\documentclass[12pt]{article}
\usepackage{amsmath,amsthm,epsfig,amsfonts,amssymb}

\theoremstyle{plain}

\newtheorem{stat}{Statement}[section]
 \newtheorem{thm}[stat]{Theorem}
 
 \newtheorem{prop}[stat]{Proposition}
 \newtheorem{lem}[stat]{Lemma}
 \newtheorem{thm*}{Theorem C.0.\!\!}
\theoremstyle{definition}
 \newtheorem{rem}[stat]{Remark}
 \newtheorem{example}[stat]{Example}

\newenvironment{pf}{\vskip 0.3cm \noindent{\it Proof }}
{\hfill $\square$ \vskip 0.3cm}

\def \BC{\mathcal B}

\def \DC{\mathcal D}

\def \FC{\mathcal F}

\def \cL{\mathcal L}

\def \cP{\mathcal P}

\def \IR{\mathbb R}

\numberwithin{equation}{section}

\begin{document}

\begin{center}
{\bf \Large A Feynman-Kac-type formula
\vskip 12pt
 for the deterministic and stochastic wave equations
 \vskip 12pt
  and other p.d.e.'s}
\vskip 16pt

Robert C. Dalang\footnote[1]{Institut de math\'ematiques, Ecole Polytechnique F\'ed\'erale, Station 8, 1015 Lausanne,
Switzerland. robert.dalang@epfl.ch}
\footnote[2]{Partially supported by the 
Swiss National Foundation for Scientific Research.},
Carl Mueller\footnote[3]{Department of Mathematics, University of Rochester,
Rochester, NY  14627, USA. cmlr@math.rochester.edu}
\footnote[4]{Partially supported by an NSF grant.}
and Roger Tribe\footnote[5]{Department of Mathematics, University of Warwick,
CV4 7AL, UK. tribe@maths.warwick.ac.uk
\vskip 12pt

   MSC 2000 Subject Classifications. Primary: 60H15; Secondary: 60H20.
\vskip 12pt
 
Keywords and phrases. Feynman-Kac formula, wave equation, probabilistic 
representation of solutions, stochastic partial differential equations,
moment formulae.} 

\end{center}


\vskip 1in
\begin{abstract}
We establish a probabilistic representation for a wide class of linear 
deterministic p.d.e.s with potential term, including the wave equation in 
spatial dimensions 1 to 3. Our representation applies to 
the heat equation, where it is related to the classical 
Feynman-Kac formula, as well as to the telegraph and beam equations. If 
the potential is a (random) spatially homogeneous Gaussian noise, then this 
formula leads to an expression for the moments of the solution.
\end{abstract}
\vskip 1in

\vfill

\eject
\section{Introduction} \label{introduction}

The purpose of this paper is to present a form of the Feynman-Kac formula which 
applies to a wide class of linear partial differential equations with a 
potential term, and, in particular, to the wave equation in dimensions 
$d \leq 3$. In the case of the heat equation, this gives an 
expression that differs from the classical Feynman-Kac formula. As an application, we 
consider a random potential term which is a spatially homogeneous Gaussian 
random field that is white in time. In this case, our approach provides a 
probabilistic representation for all product moments of the solution, which
has already shown its usefulness (see \cite{DalMuel}).

We begin by giving an informal derivation of the representation in the special 
case of the heat equation with potential, where we can contrast it with the 
classical Feynman-Kac formula.  Consider the heat equation on $\IR^d$ with a 
deterministic potential $V(t,x)$:
\begin{eqnarray}\label{heateqwp}
  \frac{\partial u(t,x)}{\partial t} &=& \frac{1}{2} \Delta u(t,x) + V(t,x)u(t,x),  \\
  u(0,x) &=& f(x).    \nonumber
\end{eqnarray}
The classical Feynman-Kac formula for the solution 
$(u(t,x),\ t \geq 0,\ x \in \IR^d)$ (see for instance \cite{KS}) states that,
under appropriate conditions on $V$ and $f$, 
\[
  u(t,x) = E^B_x\left[f(B_t)\exp\left(\int_{0}^{t}V(t-s,B_s)ds\right)\right]
\]
where $(B_t,\ t \geq 0)$ is a Brownian motion in $\IR^d$, and $E^B_x$ is the 
expectation for Brownian motion started at $B_0=x$.  

We now heuristically derive an alternative probabilistic representation 
to (\ref{heateqwp}), which will be rigorously established as a special case
of the main result in section \ref{sec3}. We start by writing Duhamel's  
formula for the solution $u(t,x)$, using the Green's function, as follows:
\begin{equation} \label{duhamel-schrodinger}
   u(t,x) = \int_{\IR^d} p_t(x-y)f(y)dy
        + \int_{0}^{t}\int_{\IR^d} p_{t-s}(x-y)V(s,y)u(s,y)dyds,
\end{equation} 
where
\begin{equation*}
   p_t(x) = \frac{1}{(2\pi t)^{d/2}}\exp\left(-\frac{|x|^2}{2t}\right).
\end{equation*}
We use (\ref{duhamel-schrodinger}) as the start of an iteration scheme.  
Substituting this expression for $u(s,y)$ back into the right hand side 
of (\ref{duhamel-schrodinger}) suggests the following series expansion for 
$u(t,x)$:  
\begin{equation}
   u(t,x) = \sum_{m=0}^{\infty}I_m(t,x),
\label{heat-expansion}
\end{equation}
where $I_0(t,x) = \int_{\IR^d} p_t(x-y)f(y)dy$ and 
\begin{equation}
I_{m+1}(t,x) = \int_{0}^{t}\int_{\IR^d} p_{t-s}(x-y) V(s,y)I_m(s,y)dyds.
\label{heat-expansion2}
\end{equation}
We wish to write an explicit expression for $I_m(t,x)$.  To begin with, let
\begin{equation*}
w(t,x) = I_0(t,x).
\end{equation*}
For convenience, let $s_{n+1}=t$ and $y_{n+1}=x$.  Then we have
\begin{equation}\label{i-expansion}
   I_m(t,x) = \int_{T_m(t)}\int_{\IR^{md}}
   \left(\prod_{k=1}^{m}p_{s_{k+1}-s_k}(y_{k+1}-y_k)V(s_k,y_k)\right)w(s_1,y_1)
\, d\bar y \, d\bar s
\end{equation}
where
\begin{equation*}
   T_m(t) = \{(s_1,\dots,s_m):  0\leq s_1\leq\dots\leq s_m\leq t\},
\end{equation*}
$d\bar y = dy_1\cdots dy_m$ and $d\bar s = ds_1\cdots ds_m$.
An alternative derivation of the series expansion (\ref{heat-expansion}) 
for $u(t,x)$ starts by expanding the exponential in the classical 
Feynman-Kac formula as a Taylor series and confirming that the terms 
correspond to the expansion (\ref{heat-expansion}). However, we will make use
of the iterative formula (\ref{heat-expansion2}) later on. 

A basic observation is that domain of integration $T_m(t)$ has volume 
$t^m/m!$, which, except for a missing exponential factor, 
is a Poisson probability. If $N(t)$ is a rate one Poisson process, then  
$P[N(t)=m]=t^ne^{-t}/m!$. Let $\tau_1 < \tau_2 < \cdots$ be 
the times of the successive jumps of the Poisson process, and let $\tau_0=0$.  
It is well known that if we condition on $N_t=m$, then the vector 
$(\tau_1,\dots,\tau_m)$ is uniformly distributed over the simplex $T_m(t)$.  
The time reversed sequence $t-\tau_m,\dots,t-\tau_1$ is also uniformly  
distributed on $T_m(t)$.  Therefore, setting $s_k=t-\tau_{m+1-k}$ and 
replacing $y_k$ by $y_{m+1-k}$, so that $y_0=x$, we may rewrite the expression 
(\ref{i-expansion}) for $I_m(t,x)$ as
\[
e^t \, E^{N} \left[\int_{\IR^{md}}
\left(\prod_{k=1}^{m}p_{\tau_{k}-\tau_{k-1}}(y_{k}-y_{k-1})
V(t-\tau_k,y_k)\right)
w(t-\tau_m,y_m) d\bar y \,\, 1_{\{N(t)=m\}} \, \right].
\]
where $E^{N}$ is the expectation with respect to the Poisson process.  But we 
can also exploit the fact that $p_{t}(x)$ is the probability density for the 
increments of a $d$-dimensional Brownian motion. Thus 
\begin{eqnarray*}
\lefteqn{ \int_{\IR^{md}}
\left(\prod_{k=1}^{m}p_{\tau_{k}-\tau_{k-1}}(y_{k}-y_{k-1})V(t-\tau_k,y_k)\right)
   w(t-\tau_m,y_m)d\bar y   }\hspace{2in} \\
&=& E^{B}_x\left[\left(\prod_{k=1}^{m}V(t-\tau_k,B_{\tau_k})\right)
  w(t-\tau_m,B_{\tau_m})\right],   
\end{eqnarray*}
where $E^{B}_x$ denotes the expectation with respect to Brownian motion started 
at $x$, and therefore,
\begin{equation*}
I_m(t,x) = e^t E^{B}_xE^{N}\left[
\left(\prod_{k=1}^{m}V(t-\tau_k,B_{\tau_k})\right)w(t-\tau_m,B_{\tau_m})
\,\, 1_{\{N(t)=m\}} \, \right].
\end{equation*}
Summing over $m$, we get
\begin{equation}
\label{main-expression}
u(t,x) = e^t E^{B}_xE^{N}\left[
\left(\prod_{k=1}^{N_t}V(t-\tau_k,B_{\tau_k})\right)
w \left(t-\tau_{N(t)},B_{\tau_{N(t)}} \right) \right].
\end{equation}

The representation (\ref{main-expression}), unlike the classical Feynman-Kac formula, 
does not use the entire Brownian path but only the values 
at a finite (random) set of times.
This allows us to extend this type of representation to equations where the 
differential operator is not the infinitesimal generator of a Markov process.
All we will require is a Poisson process and an independent 
stochastic process whose one dimensional marginals give the Green's function
for the differential operator. In particular, we will treat the case of the 
wave equation with potential in dimensions $d \leq 3$.

The outline of this paper is as follows. In Section \ref{series}, we describe 
the class of equations that we will consider and establish a series 
representation in the case of a bounded potential. In Section \ref{sec3}, we 
establish our Feynman-Kac-type formula analogous to (\ref{main-expression}),
where the Brownian motion will be replaced by a 
suitable spatial motion that depends on the particular equation being considered.
In Section \ref{sec4}, we give an application to the
situation where the potential is a Gaussian random field whose covariance is 
formally given by
\begin{equation*}
  E\left[\dot F(t,x) \dot F(s,y)\right] = \delta_0(t-s) f(x-y).
\end{equation*}
In this equation, $\delta(\cdot)$ denotes the Dirac delta function,  
$f: \IR^d \to \IR$ is continuous on $\IR^d \setminus \{0\}$ and the right-hand 
side is such that $f(x-y)$ is indeed a covariance. 
This type of covariance is widely used in the literature, including for instance 
\cite{dapratoz,D99,dalangfrangos,mueller}. 
In this case, unlike the classical Feynman-Kac formula, the noise is too 
rough for the probabilistic representation of the solution 
to make sense. Instead, we
establish in this section a formula for the second moment of the solution. 
Section \ref{sec5} contains the extension to $n$-fold product moments. This 
formula makes use of a Poisson random measure combined with a spatial motion. 
The first two named authors have made use of this formula \cite{DalMuel} to 
establish intermittency properties of the solution to the wave equation with 
potential.

We end by making a few comments on related literature. 
Probabilistic representations of the solution to deterministic 
p.d.e.'s abound. The closest related work seems to be 
results on random evolutions, surveyed for example in
Hersch \cite{Hersch} and Pinsky \cite{pinsky}. These references
give probabilistic representations for some hyperbolic equations,
including the Poisson representation for the damped wave equation in one 
spatial dimension (also known as the telegraph equation) as developed by Marc Kac
\cite{Kac1,Kac2}. Related also is the use of random flight models
for the Boltzmann equation, as described, for example, 
in \cite{pinsky,fmeleard}. 
We cannot quite find our approach represented in this
literature. 
The use of Poisson probabilities is also implicit in other 
works: Albeverio and Hoegh-Krohn \cite{albevario1} and 
Albeverio, Blanchard, Coombe,
Hoegh-Krohn, and Sirugue \cite{albevario2} have used the idea that the
multiple integrals involved in the expansion of Feynman
integrals are related to Poisson probabilities (see also \cite{wolpert}). All these works 
display the usefulness of probabilistic representations in studying 
problems of asymptotics, homogenization and perturbation theory for the
deterministic p.d.e. For parabolic equations with random potentials, the 
classical Feynman-Kac formula has been a key tool, for example in the 
parabolic Anderson problem (see Carmona and Molchanov \cite{carmmolch})
and in random waves (see Oksendal, Vage and Zhao \cite{oksendall}). 
We hope that our representation may play a similar role for other equations with 
a random potential.
%
\section{Series representation for bounded potential} \label{series}
Our probabilistic representation will be for the integral equation
\begin{equation} \label{rd0}
u(t, x) = w(t, x) + \int_0^t ds \int_{\IR^d} S(s, dy) V(t-s, x-y) u(t-s, x-y).
\end{equation}
In this section, the key assumption is the following:

\vspace{.1in}

\noindent
{\bf Assumption A.} For each $t \geq 0$, $S(t, dy)$ is a signed measure on
$\IR^d$ satisfying 
\begin{equation*}
  \sup_{t \in [0,T]} |S(t, \IR^d)| < \infty \quad \mbox{for all $T>0$,}
\end{equation*}
where $\vert S(t, \IR^d)\vert$ denotes the total variation.

\vspace{.1in}

\noindent
It is well known that a large class of linear partial differential equations of the form 
\[ Lu(t,x) = V(t,x) u(t,x)\]
 can be recast, using their Green's functions, into 
this integral form. We briefly recall some illustrative examples that 
we consider later. 
  
\begin{example}\label{example1}
\begin{itemize}
\item[(a)] {\em The heat equation on $\IR^d$.} Take $L = \frac{\partial}{\partial t} - \frac{1}{2} \Delta$
and $S(t, dy) = p_t(y) dy$. Then for a suitable initial condition $ u(0,x) = f_0(x)$,
 the Green's function 
representation of the heat equation leads to the integral equation (\ref{rd0})
with
\begin{equation*}
  w(t, x) = \int_{\IR^d} f_0(x-y) p_t(y) dy.
\end{equation*}
\item[(b)] {\em The wave equation on $\IR^d$ for $d \leq 3$.} 
Take $ L = \frac{\partial^2}{\partial t^2} - \Delta$ and 
\begin{equation*}
   S(t,dy) = \left\{\begin{array}{ll}
    \frac{1}{2} 1_{\{\vert y \vert < t\}} dy & \mbox{if } d=1,\\[8pt]
  \frac{1}{2\pi \sqrt{t^2-\vert y \vert^2}} 1_{\{\vert y \vert < t\}} dy &  \mbox{if }  d=2,\\[8pt]
  \frac{\sigma^{(2)}_t(dy)}{4 \pi t} & \mbox{if } d=3,
  \end{array}\right.
\end{equation*}
where $\sigma_t^{(2)}$ denotes the surface area on $\partial B(0,t)$ 
(the boundary of the ball centered at 0 with radius $t$). For all three values of $d$, 
$S(t, \IR^d) = t$. The initial conditions are of the form 
$ u(0,x) = f_0(x)$ and $ \frac{\partial}{\partial t} u(0,x) = f_1(x)$
for given $f_0, f_1 : \IR^d \to \IR$. In this case, letting $\ast$ 
denote convolution, 
\begin{equation*}
  w(t,x) = \frac{\partial}{\partial t} (S(t) \ast f_0)(x) + (S(t) \ast f_1) (x).
\end{equation*}
For $d \geq 4$, the fundamental solution of the wave equation is not a signed 
measure and so assumption A will not hold.
\item[(c)] {\em The wave equation with damping.}  Take
$ L = \frac{\partial^2}{\partial t^2}  + 2a \frac{\partial}{\partial t}- \Delta$ on
$\IR^d$. This also falls into the considered class when $d \leq 3$. Then 
\begin{equation*}
  S(t, dy) = \left\{
  \begin{array}{ll}
  \frac{e^{-at}}{2} I_0 (\vert a \vert \sqrt{t^2-y^2}) 1_{\{\vert y \vert < t \}} dy &  \mbox{if }  d=1,\\
  \frac{e^{-at}}{2 \pi} \frac{\cosh (\vert a \vert \sqrt{t^2-\vert y \vert^2})}{\sqrt{t^2-\vert y \vert^2}}
 1_{\{\vert y \vert <t\}}dy &  \mbox{if } d= 2,\\
  \frac{e^{-at}}{4 \pi} \left( \frac{\sigma^{(2)}_t (dy)}{t} + \vert a \vert 
\frac{I_1(\vert a \vert \sqrt{t^2-\vert y \vert^2})}{\sqrt{t^2- \vert y \vert^2}}
 1_{\{\vert y \vert < t \}} dy\right) & \mbox{if } d = 3.
  \end{array}
  \right.
\end{equation*}
In these formulas, given for instance in \cite{leveque} and \cite{dleveque}, 
$I_0$ and $I_1$ are modified Bessel functions of the first kind and of orders 
$0$ and $1$, respectively. In these three dimensions, $S(t, dy)$ is a non-negative measure.
\item[(d)] {\em The beam equation.} In dimension $d=1$, this is given by
$ L = \frac{\partial}{\partial t} 
 +  \frac{\partial^4}{\partial x^4}$ on $\IR$. 
Then $S(t,dy)=q_t(y)dy$, where $q_t(y)$ has Fourier transform $\exp(-\vert\xi\vert^4 t)$
for $t>0$. The smoothness and integrability of $q_t$, and hence assumption A,
can be deduced from the Fourier transform (for example 
$| x^2 q_t(x) | \leq C \| \partial^2 \hat{q}/\partial \xi^2 \|_1$). 
\end{itemize}
\end{example}

We now give a series representation for the solution $u(t,x)$ of (\ref{rd0}).
 
\begin{prop} Let $S(t, dy)$ be a signed measure satisfying assumption A.
Suppose that $V(t, y)$ and $w(t,x)$ are bounded measurable functions
on $[0,T] \times \IR^d$.
Define $H_0(t,x) = w(t,x),$ and, for $m \geq 0$, 
\begin{equation}\label{rd1}
H_{m+1} (t,x) = \int_0^t ds \int_{\IR^d} S(s, dy) V(t-s, x-y) H_m(t-s, x-y).
\end{equation}
Then  the integral equation (\ref{rd0}) has a unique solution 
satisfying $\sup_{t \leq T,\, x \in \IR^d} E[|u(t,x)|^2] < \infty$ given by
\begin{equation}\label{rda0}
   u(t,x) = \sum_{m=0}^\infty H_m(t, x)
\end{equation}
(the series converges uniformly on $[0,T]\times \IR^d$).
\label{prop1}
\end{prop}

\proof We first check the convergence of the series (\ref{rda0}). Set
$M_m(s) = \sup_{z \in \IR^d} \vert H_m(s,z)\vert$. Then 
\begin{eqnarray*}
M_{m+1}(t) &\leq& \sup_{z \in \IR^d} \int_0^t ds \int_{\IR^d} 
\left| S(s, dy) \right| \, \sup_{r,z} \vert V(r,z) \vert \, \sup_z \vert H_m(t-s, x-z)\vert\\
&\leq& C(S,V) \int_0^t ds\, M_m(s).
\end{eqnarray*}
A simple induction argument shows that $M_m(s) < \infty$ for all $m,s$.
Gronwall's lemma (see e.g. \cite[Remark (6)]{D99}) now implies that
$\sum^\infty_{m=0} H_m(t,x)$ converges uniformly on $[0, T] \times \IR^d.$ 

Another Gronwall argument shows the uniqueness of solutions to (\ref{rd0}).
So it suffices to check that $\sum^\infty_{m=0} H_m(t,x)$ satisfies 
(\ref{rd0}). This is the case, since 
\begin{eqnarray*}
   && w(t,x) + \int_0^t ds \int_{\IR^d} S(s, dy) V(t-s, x-y) \sum_{m=0}^\infty H_m(t-s, x-y)\\
&&\qquad = I_0 (t,x) + \sum^\infty_{m=0} \int_0^t ds \int_{\IR^d} S(s, dy) V(t-s, x-y) H_m(t-s, x-y)\\
&&\qquad = \sum^\infty_{m=0} H_m(t,x).
\end{eqnarray*}
Fubini's theorem, used for the first equality, applies by uniform 
convergence of the series and assumption A.  
\hfill $\Box$
\bigskip

\section{Probabilistic representation} \label{sec3}
For the probabilistic representation, we use the following additional assumption on 
the kernel $S(t, dy)$ used in the integral equation
(\ref{rd0}).

\vspace{.1in}

\noindent
{\bf Assumption B.} 
There exists a jointly measurable process 
$(\tilde{X}_t,\ t > 0)$ such that 
for each $t > 0$  
\begin{equation*}
P \left\{-\tilde{X}_t\in dx \right\} = \frac{|S(t, dx)|}{|S(t, \IR^d)|}.
\end{equation*}
(In the case where $S(t, A) = S(t, -A),$ the minus sign in front of $\tilde{X}_t$ is not needed).

\vspace{.1in}

\begin{example} 
\begin{itemize}
\item[(a)] {\em The heat equation.} In this case, 
one can take $\tilde{X}_t = \sqrt{t} \ X_0,$ where $X_0$ is a standard 
$N(0, I_d)$ random vector in $\IR^d$. An alternative possibility is to let 
$(\tilde{X}_t)$ be a standard Brownian motion in $\IR^d$. 
\item[(b)] {\em The wave equation.} In the three dimensional case, 
one can take $\tilde{X}_t = t \, \Theta_0,$ where $\Theta_0$ is chosen according to 
the uniform probability measure on $\partial B(0,1).$ 
The one and two dimensional cases can be handled in a similar way.
\label{rdexample2.1}
\item[(c)] {\em The damped wave equation.} Kac \cite{Kac1,Kac2} pointed out a neat representation for 
solutions to the damped wave equation in dimension 1. Let $(N_a(t))$ be a rate $a$ Poisson process
and define $\tau_t = \int^t_0 (-1)^{N_a(s)} ds$. If $w(t,x)$ solves the undamped wave equation
$ \frac{\partial^2 w}{\partial t^2}  - \Delta w=0$ for $t \in \IR$, $x \in \IR^d$,
then $u(t,x) = E [w(\tau_t,x)]$ solves the damped wave equation
$\frac{\partial^2 u}{\partial t^2}  + 2a \frac{\partial u}{\partial t}- \Delta u=0$,
and with the same initial conditions. Using this, one finds that the kernel $S_a(t,dy)$
for the damped equation can be written as $S_a(t,dy) = 
E[ S(\tau_t,dy)1_{\{\tau_t >0\}} + S(-\tau_t,dy) 1_{\{\tau_t<0\}}]$, where $S(t,dy)$ is the kernel for the wave equation. Now
we satisfy assumption B by setting $\tilde{X}_t= |\tau_t| \, \Theta_0$,
where $\Theta_0$ is a uniform random variable on $[-1,1]$, independent of $N_a(t)$. 
\item[(d)] {\em The beam equation.} As in example (a), we can use scaling to set $\tilde{X}_t
= t^{1/4} X_0$, where $X_0$ is chosen to have distribution $|S(1,dy)|$.
\end{itemize}
\end{example}
\vskip 12pt

Let $\tilde{X}^{(i)} = (\tilde{X}_t^{(i)},\ t \geq 0),\ i \geq 1,$ be 
i.i.d.~copies of $(\tilde{X}_t,\ t \geq 0)$, and let $(N(t),\ t \geq 0)$ be a rate 
one Poisson process independent of the $(\tilde{X}^{(i)}).$ Let 
$0 < \tau_1 < \tau_2 < \cdots$ be the jump times of $(N(t))$ and set 
$\tau_0 \equiv 0.$ 
Define a process $X = (X_t,\ t \geq 0)$ as follows :
\begin{equation*}
   X_t = X_0 + \tilde{X}^{(1)}_t \mbox{ for } 0 \leq t \leq \tau_1,
\end{equation*}
and for $i \geq 1,$
\begin{equation*}
   X_t = X_{\tau_i} + \tilde{X}^{(i+1)}_{t-\tau_i}, \mbox{ for } \tau_i < t \leq \tau_{i+1}.
\end{equation*}

\begin{figure}[h]
\begin{picture}(400,180)(0,0)
\put(120,20){\includegraphics[scale=1,angle=0]{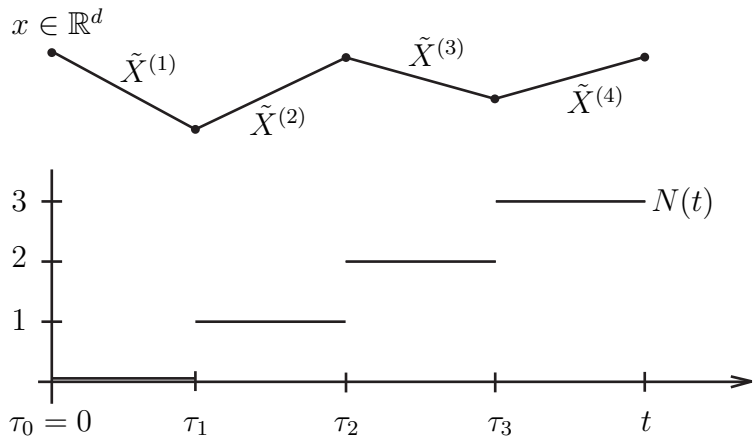}}
\put(109,10){$\tau_0 = 0$}
\put(175,10){$\tau_1$} \put(232,10){$\tau_2$}
\put(290,10){$\tau_3$}
\put(348,10){$t$}
\put(110,160){$x \in \IR^d$}
\put(151,142){$\tilde X^{(1)}$}
\put(200,122){$\tilde X^{(2)}$}
\put(260,148){$\tilde X^{(3)}$}
\put(320,130){$\tilde X^{(4)}$}
\put(110,48){1}
\put(110,72){2}
\put(110,93){3}
\put(352,93){$N(t)$}
\end{picture}
\caption{A sample path of the process $X$ and of the Poisson process $(N(t))$. \label{fig1} }
\end{figure}

   We use $P_x$ to denote a probability under which, in addition,
$X_0=x$ with probability one. Informally, the process $X$ follows $\tilde{X}^{(1)}$ during the interval 
$[0, \tau_1]$, then follows $\tilde{X}^{(2)}$ started at $X_{\tau_1}$ during 
$[\tau_1, \tau_2]$, then $\tilde{X}^{(3)}$ started at $X_{\tau_2}$ during 
$[\tau_2, \tau_3]$, etc. See Figure \ref{fig1} for an illustration.

\begin{thm} \label{maintheorem}
Suppose that the kernel $S(t,dy)$ is a non-negative measure
satisfying assumptions A and B. Suppose $w(t,x)$ is bounded and measurable for
$t \leq T$, $x \in \IR^d$. Then 
$(u(t,x),\ t \leq T, x \in \IR^d)$ defined by  
\begin{equation}
\label{expression-for-u}
   u(t,x) = e^t E_x \left[w \left(t-\tau_{N(t)}, X_{\tau_{N(t)}} \right) 
\prod_{i=1}^{N(t)} \Big[ S(\tau_i-\tau_{i-1}, \IR^d) 
V(t-\tau_i, X_{\tau_i})\Big]\right]
\end{equation}
(where, on $\{N(t) = 0\}$, the product is defined to take the value $1$)
is the solution of (\ref{rd0}).
\label{thm2}
\end{thm}

\proof For $t \geq 0$ and $m \geq 0$, let
\begin{equation*}
   Y(m,t) = e^t 1_{\{N(t) = m\}} w(t-\tau_m, X_{\tau_m}) \prod^m_{i=1}
     \Big[S(\tau_i-\tau_{i-1}, \IR^d) V(t-\tau_i, X_{\tau_i})\Big].
\end{equation*}
Then $ u(t,x) = \sum^\infty_{m=0} E_x [Y(m,t)]$. In order to show that $u(t,x)$ is 
the solution of (\ref{rd0}), it suffices by Proposition \ref{prop1} to show that $H_m(t,x) = E_x[Y(m,t)]$ 
for all $m,t,x$. We prove this by induction on $m$. For $m=0$, 
\begin{eqnarray*}
  E_x[Y(m,t)] & = & E_x \left[e^t 1_{\{N(t) = 0\}} w(t, X_0)\right]\\
   &=& e^t w(t,x) P_x\{N(t) = 0\} = w(t,x) = H_0(t,x).
\end{eqnarray*}

Now fix $m \geq 1$ and suppose by induction that $H_{m-1}(t,x) = E_x[Y(m-1, t)]$, 
for all $(t,x)$. Set $\FC_1 = \sigma\{ X_{\tau_1}, \, \tau_1\}$. Then
\begin{eqnarray*}
  && E_x[Y(m,t)] \\
&&\qquad =  E_x \bigg[S(\tau_1, \IR^d) V(t-\tau_1, X_{\tau_1}) 1_{\{\tau_1 \leq t\}} e^{\tau_1} \\
   &&\qquad \qquad \qquad \times \, E_x \Big[e^{t-\tau_1} 1_{\{N(t)-N(\tau_1)=m-1\}} 
w ((t-\tau_1)-(\tau_m-\tau_1), X_{\tau_m}) \\
  && \qquad \qquad \qquad \qquad \times \, \prod^m_{i=2} \big\{S((\tau_i-\tau_1)-(\tau_{i-1}-\tau_1), \IR^d)
 V((t-\tau_1)-(\tau_i-\tau_1), X_{\tau_i})\big\} \bigg\vert \FC_1 \Big]\bigg].
\end{eqnarray*}
Note that, for $i \geq 1$, 
\[
X_{\tau_i} = X_{\tau_1} + \sum_{j=1}^{i-1} \tilde{X}^{(j+1)}_{\tau_{j+1}-\tau_j},
\]
and the conditional expectation can be expressed using only the increments $\tau_i-\tau_1$ for
$i \geq 1$.  Using the strong Markov property of $(N(t))$ at time $\tau_1$ and the 
independence of the families $X^{(i)}_t$, we deduce that 
\begin{eqnarray*}
E_x[Y(m,t)]
  &=& E_x \bigg[S(\tau_1, \IR^d) V(t-\tau_1, X_{\tau_1}) e^{\tau_1} 1_{\{\tau_1 
  \leq t\}} Y_{m-1}( t-\tau_1, X_{\tau_1})\bigg]  \\
   &=&  \int_0^t ds \ e^{-s} S(s, \IR^d) e^s \int_{\IR^d} \frac{S(s, dy)}{S(s, \IR^d)} 
V(t-s, x-y) Y_{m-1} (t-s, x-y)\\
&=& \int_0^t ds \int_{\IR^d} S(s, dy) V(t-s, x-y) H_{m-1}(t-s, x-y)\\
&=& H_m(t,x),
\end{eqnarray*}
by the induction hypothesis and (\ref{rd1}). This completes the proof. 
\hfill  $\Box$
\bigskip

We have presented the simplest setting of the probabilistic representation, 
sufficient to treat our interest in the wave equation in dimensions $d \leq 3$ 
and the random potentials in the subsequent sections. However various 
extensions and variations of this representation are possible. We give a brief 
description here, leaving the details to the interested reader. 
\begin{enumerate}
\item For a signed kernel $S(t,dy)$, we need to modify somewhat the representation.
Write $S(t,dy) = S_{+}(t,dy) - S_{-}(t,dy)$ for the Hahn-Jordan decomposition into a difference of non-negative measures. Choose, if possible, subsets $A(t) \subseteq \IR^d$
so that $S_+(t,A(t)) = S_-(t,A^c(t)) =0$ and $(x,t) \to 1_{A(t)}(x)$ is measurable.
(Note that this is certainly possible when $S(t,dy) =q_t(y)dy$ for continuous 
$(q_t(y),\ t>0,\ y \in \IR^d)$.) Let $C_t$ be a counter defined by
\[
C_t = \sum_{i=1}^{\infty} 1_{\{ \tau_i \leq t\}} \, 1_{\{ X_{\tau_{i-1}}- X_{\tau_i} \in 
A(\tau_i- \tau_{i-1}) \}}.
\]
Then the argument above leads to the representation
\[
u(t,x) = e^t E_x \left[ w \left(t-\tau_{N(t)}, X_{\tau_{N(t)}} \right) (-1)^{C_t} 
\prod_{i=1}^{N(t)} \Big[ \left| S(\tau_i-\tau_{i-1}, \cdot) \right| 
V(t-\tau_i, X_{\tau_i})\Big]\right],
\]
where $|S(t,\cdot)|$ denotes the total variation of the measure $S(t,dy)$.  
This representation then covers the case of the beam equation in all dimensions $d\geq 1$.
\item If, instead of being real-valued, 
$u(t,x) = (u_1(t,x), \dots,u_n(t,x)) \in \IR^n$, and $V(t,x)$ is an 
$n \times n$ matrix, so that (\ref{rd0}) is in fact a system of p.d.e.'s, then 
the formula in Theorem \ref{thm2} still holds, provided the matrix product 
in (\ref{expression-for-u}) is ordered according to increasing values of $i$. 
\item We have treated for simplicity the case of spatially homogeneous equations on 
$\IR^d$. However, in principle, suitable changes should allow representations for 
inhomogeneous equations, or equations in domains with suitable boundary 
conditions.
\item For any $\lambda>0$, one can replace the potential $V$ by 
$\lambda^{-1} V$ and use a Poisson process of rate $\lambda$ to obtain an 
alternative representation.  For example, rewriting the heat equation 
(\ref{heateqwp}) as
\[
 \frac{\partial u}{\partial t} = \frac12 \Delta u  
    + \lambda \left[\frac{V}{\lambda} \right]u,
\]
we would get the representation
\begin{equation}
\label{lambda-rep}
u(t,x) = e^{\lambda t}
E_x \left[ w \left(t-\tau_{N(t)}, X_{\tau_{N(t)}} \right) (-1)^{C_t} 
\prod_{i=1}^{N(t)} \Big[ \left| S(\tau_i-\tau_{i-1}, \cdot) \right| 
\lambda^{-1} V(t-\tau_i, X_{\tau_i})\Big]\right]
\end{equation}
where $X_t$ starts over at the times of a rate $\lambda$ Poisson process.  For 
large $\lambda$, these representations, when using a Markovian $X$, become 
close to the classical Feynman-Kac formula.  For example, we can further 
rewrite (\ref{heateqwp}) as
\begin{equation}
\label{new-heat}
 \frac{\partial u}{\partial t} = \frac12 \Delta u  + \lambda \left[1+ \frac{V}{\lambda} \right]u - \lambda u.  
\end{equation}
Due to the term $-\lambda u$ in (\ref{new-heat}), the Green's function $e^{-\lambda t} p_t(y)$ of $Lu = \frac{\partial}{\partial t} - \frac{1}{2} \Delta u + \lambda u$ gives rise to a factor 
$e^{-\lambda t}$ inside the expectation in (\ref{lambda-rep}), which cancels the factor
$e^{\lambda t}$ which is outside of the expectation. Regarding $[1+(V/\lambda)]$ as our potential term, we find
\[
u(t,x) = E_x^B \left[ w(t-\tau_{N_{\lambda}(t)},B_{\tau_{N_{\lambda}(t)}}) \exp \left(
\sum_{m=1}^{N_{\lambda}(t)} \log \left( 1+  \lambda^{-1} V(t-\tau_m,B_{\tau_m}) 
\right) \right) \right]. 
\]
Now letting $\lambda \to \infty$ and using $\log(1+x) \approx x$, 
the integrand involves a Riemann sum approximation to the integral in the 
classical Feynman-Kac formula.
\end{enumerate}

\section{Second moments for random potentials} \label{sec4} 
\subsection{The random potentials} \label{sec4.1}
We are now going to consider a class of linear 
equations driven by spatially homogeneous Gaussian noise $\dot F(t,x),$ whose 
covariance is formally given by 
\begin{equation*}
  E\left[\dot F(t,x) \dot F(s,y)\right] = \delta_0(t-s) f(x-y).
\end{equation*}
In this equation, $\delta(\cdot)$ denotes the Dirac delta function, and 
$f: \IR^d \to \IR$ is continuous on $\IR^d \setminus \{0\}$. More precisely, 
let $\DC(\IR^{d+1})$ be the space of Schwartz test functions (see 
\cite{schwartz}). On a given probability space, we define a Gaussian process 
$F = (F(\varphi),\ \varphi \in \DC(\IR^{d+1}))$ with mean zero and covariance 
functional 
  \[
  E\big[F(\varphi)F(\psi)\big] = \int_{\IR_+} dt \int_{\IR^d} dx \int_{\IR^d} dy \, \varphi(t,x) f(x-y) \psi(t,y).
 \]
Since this is a covariance, it is well-known \cite[Schwartz, Chap. VII, 
Th\'eor\`eme XVII]{schwartz} that $f$ must be symmetric and be the Fourier 
transform of a non-negative tempered measure $\mu$ on $\IR^d$, termed the 
spectral measure : $f = \FC \mu$. In this case, $F$ extends to a worthy 
martingale measure $M = (M_t(B),\ t \geq 0,\ B \in \BC_b(\IR^d))$ in the sense 
of \cite{walsh}, with covariation measure $Q$ defined by 
 \begin{equation*}
 Q([0,t] \times A \times B) = \langle M(A), M(B)\rangle_t = t \int_{\IR^d} dx \int_{\IR^d} dy \  1_A(x) f(x-y) 1_B(y),
 \end{equation*}
and dominating measure $K = Q$ (see \cite{dalangfrangos,D99}). By 
construction, $t \mapsto M_t(B)$ is a continuous martingale and 
\begin{equation*}
 F(\varphi) = \int_{\IR_+ \times \IR^d} \varphi(t,x) M(dt, dx),
\end{equation*}
where the stochastic integral is as defined in \cite{walsh}.
\vskip 12pt

\noindent{\bf Assumption C.} For each $t>0$, $S(t,dy)$ is a non-negative measure and takes 
values in the space of distributions with rapid decrease \cite[Chap.VIII, \S 5]{schwartz}. 
Moreover, it satisfies
\begin{equation} \label{rd3.2}
   \int^T_0 ds \int_{\IR^d} \mu(d \xi)\, \vert \FC S(s,\cdot)(\xi)\vert^2 < \infty
\end{equation}
and
\begin{equation*}
   \lim_{h\downarrow 0} \int_0^T dt \int_{\IR^d} \mu(d\xi)\, \sup_{t<r<t+h} 
\vert \FC S(r)(\xi) - \FC S(t)(\xi)\vert^2 = 0.
\end{equation*}
\vskip 12pt

\noindent
Consider the stochastic integral equation 
\begin{equation}\label{rd3.1}
   u(t,x) = w(t,x) + \int_0^t \int_{\IR^d} S(t-s, x-y) u(s,y) F(ds, dy),
\end{equation}
where $w(t,x)$ is a random field satisfying appropriate conditions (see below).

   Our motivation is the case where $S(t,dy)$ is the Green's function for 
a partial differential operator $L$, and the study of the stochastic p.d.e.
$ L u = u \ \dot{F}$, with stationary initial conditions independent of $\dot F$. This s.p.d.e.~can be recast into this integral form 
with $w(t,x)$ being the solution of $Lw = 0$ with the same initial conditions as $u(t,x)$. In this context, $(w(t,\cdot),\ M_t(\cdot))$ is stationary in $x$, or, more precisely, has property (S) of Dalang \cite[Definition 5.1]{D99}.

   The stochastic integral in (\ref{rd3.1}) needs defining. If $S(s,y)$ is a smooth function, as in the heat equation, then we can use the stochastic integral with respect to a worthy martingale measure introduced in \cite{walsh}. In this case, (\ref{rd3.1}) has a unique solution provided $w(t,x)$ is a predictable process such that $\sup_{t \leq T,\ x \in \IR^d} E[w^2(t,x)] < \infty$.
   
   If $S(s,\cdot)$ is a singular measure, as in the case of the 3-dimensional wave equation that we are particularly interested in, then we use the integral introduced in Dalang \cite{D99}. We briefly describe his construction, that uses an approximation to the identity. Choose $\psi \in C_0^\infty(\IR^d)$ with $\psi \geq 0$, the support of $\psi$ is 
contained in the unit ball of $\IR^d$ and $\int_{\IR^d} \psi(x) dx =1$. For 
$\ell \geq 1$, set $\psi_\ell(x) = \ell^d \psi(\ell x)$, so that $\psi_\ell \to \delta_0$ as 
$\ell \to \infty$. The stochastic integral in (\ref{rd3.1})
is the $L^2$-limit of the usual stochastic integrals
\[
\int_0^t \int_{\IR^d} S_\ell(t-s, x-y) u(s,y) F(ds, dy).
\] 
where $S_\ell(t,x)$ is the convolution $\int S(t, dy) \psi_\ell(x-y)$. While studying the s.p.d.e.~$Lu=u \dot F$ as above, this convergence was established in \cite{D99}, and the same arguments show existence and uniqueness of a solution to (\ref{rd3.1}) provided $w(t,x)$ has the property (S) of \cite[Definition 5.1]{D99} and $\sup_{t \leq T} E[w^2(t,0)] < \infty$ for all $T>0$.
Details for this can be found in \cite[Section 5]{D99}. Assumption C (in particular, the fact that $S(t,dy)$ is non-negative) is also used in the definition of the stochastic integral. In the terminology of \cite{dleveque}, $u(t,x)$ is a {\em random field solution} of (\ref{rd3.1}), that is defined for every $t$ and $x$ (as opposed to a {\em function-valued solution}, defined only for all $t$ and almost all $x$, that would not be adequate for our purposes).

   In fact, it is shown in \cite{D99} that (\ref{rd3.2}) is even a necessary condition for (\ref{rd3.1}) to have a solution 
satisfying $\sup_{t \leq T,\, x \in \IR^d} E[u^2(t,x)] < \infty$.

   In the cases of the heat and wave
equations, \cite{D99} gives equivalent conditions to (\ref{rd3.2})
involving only $\mu$ or the function $f$ in the covariance structure.

\subsection{The series representation}\label{sec4.2}

In this subsection, we work under assumptions A and C. We assume that $w(t,x)$ has the properties indicated in subsection \ref{sec4.1} that ensure that the stochastic integral in (\ref{rd3.1}) is well-defined and ensure existence and uniqueness of a random field solution to this integral equation. 

   We shall show that there is a series representation for the solution $u(t,x)$ of 
(\ref{rd3.1}), analogous to (\ref{rda0}), but with the deterministic integral 
replaced by the stochastic integral as in (\ref{rd3.1}).    
Define $I_0(t,x) = w(t,x),$ and, for $m \geq 0$,
\begin{equation}  \label{rd3.3a}
   I_{m+1} (t, x) = \int_0^t \int_{\IR^d} S(t-s, x-y) I_m(s, y) F(ds, dy).
\end{equation}
\begin{prop} Suppose that $w(t,x)$ is bounded and measurable for 
$t \leq T, \, x \in \IR^d$. Then the series 
\begin{equation}\label{rd3.3}
   u(t,x) = \sum^\infty_{m=0} I_m (t, x)
\end{equation}
converges in $L^2$ uniformly over $(t,x) \in [0,T] \times \IR^d$
and is the unique solution to (\ref{rd3.1}). 
\label{prop3.1}
\end{prop}

\proof We first check the $L^2$-convergence of the series in (\ref{rd3.3}). 
Set 
\begin{equation*}
M_m(t) = \sup_{x \in \IR^d} E[I_m(t, x)^2].
\end{equation*}
By \cite[Theorem 2]{D99},  
\begin{equation*}
M_m(t) \leq \int_0^t ds\, M_{m-1} (s) \int_{\IR^d} \mu(d\xi)\, \vert \FC S(t-s,\cdot)(\xi)\vert^2,
\end{equation*}
 By (\ref{rd3.2}) and \cite[Lemma 15]{D99,D99b}, we conclude that
\begin{equation*}
   \sum^\infty_{m=0} M_m(s)^{1/2} < \infty,
\end{equation*}
which establishes the $L^2$-convergence of the series. Set 
$u_n(t,x) = \sum^n_{m=0} I_m(t,x)$.  
Then $u_n(t,x) \to u(t,x)$ in $L^2$, and by \cite[Theorem 2]{D99}, as 
$n\to\infty$, 
\begin{equation*}
\int^t_0 \int_{\IR^d} S(t-s, x-y) u_n(s, y) F(ds, dy) \stackrel{L^2}{\to} \int_0^t \int_{\IR^d} S(t-s, x-y) u(x, y) F(ds, dy).
\end{equation*}
Therefore,
\begin{eqnarray*}
   && w(t,x) + \int^t_0 \int_{\IR^d} S(t-s, x-y) u(s, y) F(ds, dy)\\
   && \qquad = \lim_{n \to \infty} \left(w(t,x) + \int^t_0 \int_{\IR^d} S(t-s, x-y) u_n (s,y) F(ds, dy) \right)\\
   && \qquad = \lim_{n \to \infty} \left(I_0(t,x) + \sum^n_{m=0} \int_0^t \int_{\IR^d} S(t-s, x-y) I_m(x, y) F(ds, dy)\right) \\
   && \qquad = \lim_{n \to \infty} \sum^{n+1}_{m=0} I_m (t,x)\\
   && \qquad = u(t,x),
\end{eqnarray*}
showing that $u(t,x)$ solves (\ref{rd3.1}). \hfill $\Box$
\bigskip

The successive terms in (\ref{rd3.3}) are orthogonal in $L^2$, that is
$E[ I_m(t,x) I_{m'}(s,y)]=0$ whenever $m \neq m'$. The series is therefore 
a chaos expansion for the noise $F$. The orthogonality can be checked
by induction on $m$ and $m'$, using the fact that the covariance between
$I_m$ and $I_{m'}$ reduces, as in (\ref{Imeqn}) below,
to an expression involving the covariance between $I_{m-1}$ and
$I_{m' -1}$. 

\subsection{The probabilistic representation of second moments}

In this subsection, we work under assumptions A, B, and C. We make the same assumptions on $w(t,x)$ as in subsection \ref{sec4.2}.

\begin{figure}[h]
\begin{picture}(400,180)(0,0)
\put(120,20){\includegraphics[scale=1,angle=0]{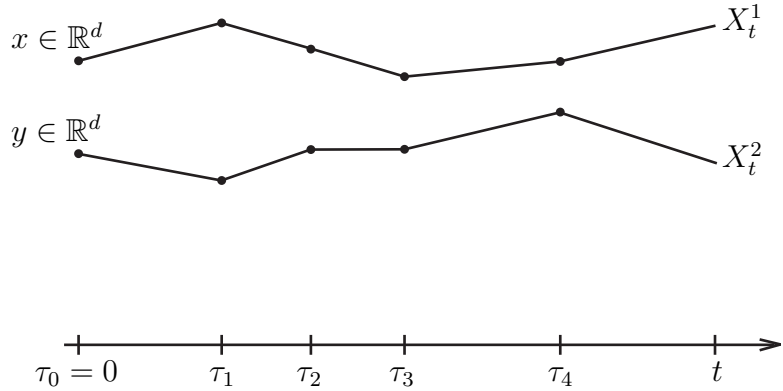}}
\put(109,10){$\tau_0 = 0$}
\put(175,10){$\tau_1$} \put(208,10){$\tau_2$}
\put(243,10){$\tau_3$}\put(303,10){$\tau_4$}
\put(365,10){$t$}
\put(100,136){$x \in \IR^d$}
\put(100,101){$y \in \IR^d$}
\put(369,143){$X^1_t$}
\put(369,91){$X^2_t$}
\end{picture}
\caption{A sample path of the processes $X^1$ and $X^2$. \label{fig2} }
\end{figure}

   Let $(N(t),\ t \geq 0)$ be a rate one Poisson process. Using two independent 
i.i.d.~families $(\tilde{X}^{(i,1)}_{\cdot}, \, i  \geq 1)$ and 
$(\tilde{X}^{(i,2)}_{\cdot}, \, i \geq 1)$, construct, as in Section
\ref{sec3}, two processes $X^1 = (X^1_t,\ t \geq 0)$ and $X^2 = (X^2_t,\ t \geq 0)$
which renew themselves at the same set of jump times
$\tau_i$ of the process $N$, and which start, under $P_{x_1,x_2}$, at $x_1$ and $x_2$ respectively.  See Figure \ref{fig2} for an illustration.

\begin{thm} Let $u(t,x)$ be the solution of (\ref{rd3.1}) given in Proposition \ref{prop3.1}. Then
\begin{eqnarray*}
E[u(t,x)u(t,y)] &=& e^{t} E_{x,y} \Big[ w\left(t-\tau_{N(t)}, X^1_{\tau_{N(t)}}\right) 
       w\left(t-\tau_{N(t)}, X^2_{\tau_{N(t)}}\right) \\
&& \qquad \qquad \times \,  \prod^{N(t)}_{i=1} \left( S(\tau_i-\tau_{i-1}, \IR^d)^2 
         f\left(X^1_{\tau_i} - X^2_{\tau_i}\right) \right) \Big].
\end{eqnarray*}
\end{thm}

\proof Observe that by Proposition \ref{prop3.1},
\begin{equation*}
   E[u(t,x) u(t,y)] = \sum^\infty_{m=0} \sum^\infty_{m'=0} E[I_m(t,x) I_{m'}(t,y)] = \sum_{m=0}^\infty E[I_m(t,x) I_m(t,y)],
\end{equation*}
using the orthogonality of the terms in the series. 
For $m \geq 1$, using the smoothed kernels $S_\ell(t,x)$ defined earlier, 
we have 
\begin{eqnarray}
&& \hspace{-.3in} E[I_m(t,x) I_m(t,y)]  \nonumber \\
& = & E\left[ \int^t_0 \int_{\IR^d} S(t-s, x-z) I_{m-1}(s,z) F(ds, dz) \int_0^t \int_{\IR^d} S(t-s, y-z) I_{m-1}(s,z) 
  F(ds,dz)\right] \nonumber \\
&=& \lim_{\ell \to \infty} E\left[ \int^t_0 \int_{\IR^d} S_\ell(t-s, x-z) I_{m-1}(s,z) F(ds, dz) 
\int_0^t \int_{\IR^d} S_\ell(t-s, y-z) I_{m-1}(s,z)   F(ds,dz)\right] \nonumber \\
&=& \lim_{\ell \to \infty} \int_0^t ds \int_{\IR^d} dz_1 \int_{\IR^d} dz_2\, 
      S_\ell(t-s, x-z_1) S_\ell(t-s, y-z_2) f(z_1-z_2) E[I_{m-1}(s, z_1) I_{m-1}(s, z_2)] \nonumber \\
&=& \int_0^t ds \int_{\IR^d} S(t-s, x-dz_1)\int_{\IR^d} S(t-s, y-dz_2) f(z_1-z_2) E[I_{m-1}(s, z_1)I_{m-1}(s, z_2)],
 \label{Imeqn}
\end{eqnarray}
where we have used the Lebesgue Differentiation Theorem \cite[Chapter 7, Exercise 2]{WZ} in the final step.
We shall now show by induction that
\begin{equation}\label{rd3.4}
E[I_m(t,x) I_m(t,y)] = J(m,t,x,y), \qquad m \geq 0, 
\end{equation}
where
\begin{eqnarray*}
J(m,t,x,y) & =&  e^{t} E_{x,y}\bigg[1_{\{N(t)=m\}} w\left(t-\tau_m, X^1_{\tau_m}\right) w\left(t-\tau_m, X^2_{\tau_m}\right)\\
&& \qquad \qquad \qquad
\times \, \prod^m_{i=1} \left\{S(\tau_i-\tau_{i-1}, \IR^d)^2 f(X^1_{\tau_i}-X^2_{\tau_i})\right\}\bigg].
\end{eqnarray*}
For $m=0,$
\begin{equation*}
J(0, t, x,y) = e^{t} w(t,x)w(t,y) P_{x,y}\{N(t) = 0\} = E[I_0(t,x) I_0 (t,y)].
\end{equation*}
We suppose now that (\ref{rd3.4}) holds for $m-1$. By the Markov property at $\tau_1$, 
arguing as in Theorem \ref{thm2}, we have, choosing 
$\FC_1 = \sigma \{ \tau_1, X^1_{\tau_1}, X^2_{\tau_1} \}$, 
\begin{eqnarray*}
 J(m,t,x,y) &=&  E_{x,y} \Big[ 1_{\{\tau_1 \leq t\}} e^{\tau_1}
       f(X^1_{\tau_1}-X^2_{\tau_1})S(\tau_1, \IR^d)^2\\
   && \qquad \qquad  \times \, e^{(t-\tau_1)}E_{x,y} [1_{\{N(t)-N(\tau_1)=m-1\}} 
     w(t-\tau_m, X^1_{\tau_m})w(t-\tau_m, X^2_{\tau_m}) \\
   && \hspace{0.9in}  \times  \, \prod^m_{i=2}[S(\tau_i-\tau_{i-1}, \IR^d)^2
    f(X^1_{\tau_i}-X^2_{\tau_i}) \vert \FC_{1} ] \Big]\\
   &=& E_{x,y} \left[1_{\{\tau_1 \leq t\}} e^{\tau_1} f(X_{\tau_1}-X^2_{\tau_1}) 
      S(\tau_1, \IR^d)^2 J(m-1, t-\tau_1, X^1_{\tau_1}, X^2_{\tau_2}) \right] \\
   &=& \displaystyle\int^t_0 ds\, \int_{\IR^d}
    S(s, x-dz_1) \int_{\IR^d} S(s, y-dz_2) f(z_1 -z_2) J(m-1, t-s, z_1, z_2).
\end{eqnarray*}
The conclusion now follows from (\ref{Imeqn}) and the induction hypothesis. 
\hfill $\Box$
\bigskip

\begin{rem} 
By multiplying the integral formulas (\ref{rd3.1}) for $u(t,x)$ and $u(t,y)$ and taking
expectations, one expects formally the integral equation 
\begin{eqnarray*}
&& \hspace{-.3in} E[u(t,x)u(t,y)] \\ 
& = & w(t,x) w(t,y) 
+ \int^t_0 \int_{\IR^d} \int_{\IR^d} S(t-s,x- dz_1) S(t-s,y- dz_2)
f(z_1- z_2) E[ u(s,z_1) u(s,z_2)].
\end{eqnarray*}
This new integral equation on $\IR^{2d}$ is of the same form as (\ref{rd0}). This 
leads to an alternative derivation, by applying Theorem \ref{maintheorem},
of the representation for second moments given above. However, we have 
used the argument above 
as it will generalize to higher moments. 
\end{rem}

\section{Moments of order $n$}\label{sec5}

In this subsection, we work again under assumptions A, B, an C. In addition to the assumptions on $w(t,x)$ made in subsection \ref{sec4.2}, we assume that 
$ \sup_{t \leq T} E[|w(t,x)|^p] < \infty$,  for all $T,p>0$, which ensures that the solutions have finite $p$-th moments.

In the case where $u(t,x)$ solves a first order equation driven by the Gaussian noise $\dot F$, 
written in the form $\partial u/\partial t = Lu + u \dot F$, then a formal calculation suggests that the $n$-th moment $m(t,x_1, \ldots,x_n) = E[u(t,x_1) \ldots u(t,x_n)] $ should satisfy
\[
\frac{\partial m}{\partial t} = 
L_{x_1,\ldots,x_n} m + \frac12 m 
\sum_{i \neq j}^n f(x_i-x_j)  
\] 
(this formula is proved for discrete space in \cite[Section II.3]{carmmolch}).
Here, $L_{x_1,\ldots,x_n}$ stands for the sum of the operator $L$ applied to each 
variable $x_i$.
The equation for $m$ is again of the same potential type considered in section
\ref{series}, and can be recast as an integral equation using a multiple product kernel 
constructed out of the kernel $S(t,dy)$ for $L$. Theorem \ref{maintheorem}
then leads to a probabilistic representation for $m$. This argument does not
seem to apply for second order equations or directly for integral equations.
However, as we shall now explain, it is possible 
to find a representation, analogous to the one for
second moments, that holds for the higher moments of the 
general integral equation (\ref{rd3.1}).

We start with an informal discussion of the representation for higher moments.
The second moments were given in terms of a pair of processes, both of which 
were renewed at the times $\tau_i$ of a single Poisson process $N$.  The situation for the 
$n$-th moment is somewhat analogous.  Instead of two processes, we will use  
$n$ processes $X^{1},\dots,X^{n}$.  For each pair of indices 
$\rho=\{\rho_1,\rho_2\}$, we create a Poisson process $N_t(\rho)$.  The renewal 
times of the the process $X^{i}$ will be the union of the Poisson times 
arising from the processes $N_t(\rho)$, such that the index $i$ is contained 
in the pair of indices $\rho$.  

\begin{figure}[h]
\begin{picture}(400,250)(0,0)
\put(120,20){\includegraphics[scale=1,angle=0]{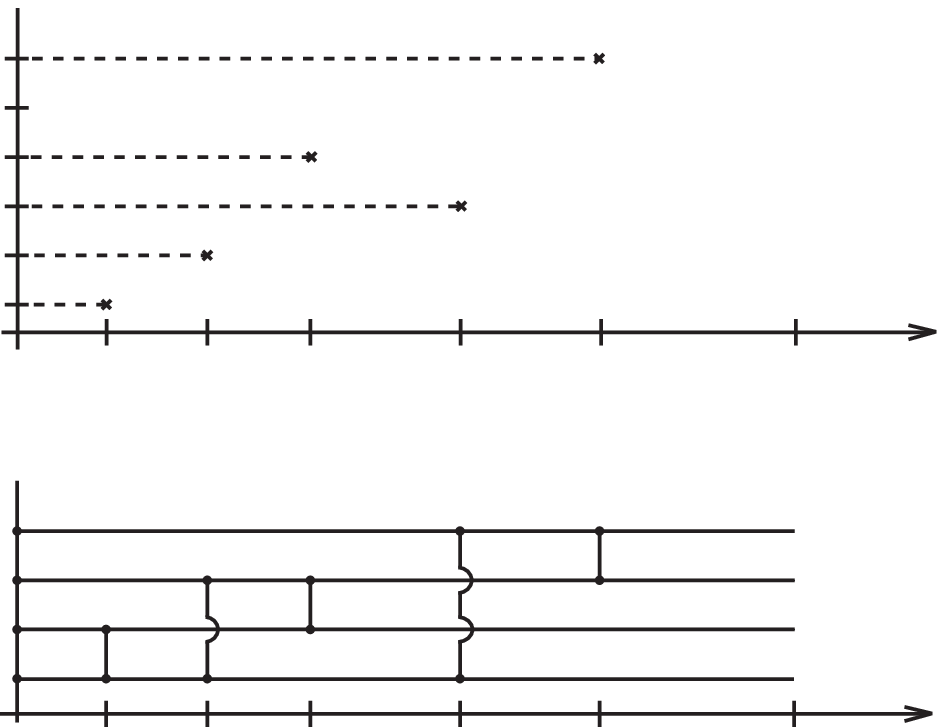}}
\put(145,10){$\sigma_1$}
\put(175,10){$\sigma_2$}
\put(205,10){$\sigma_3$}
\put(248,10){$\sigma_4$}
\put(288,10){$\sigma_5$}
\put(348,10){$t$}
\put(145,120){$\sigma_1$}
\put(175,120){$\sigma_2$}
\put(205,120){$\sigma_3$}
\put(248,120){$\sigma_4$}
\put(288,120){$\sigma_5$}
\put(348,120){$t$}
\put(115,30){1}
\put(115,45){2}
\put(115,60){3}
\put(115,75){4}
\put(90,139){$\{1,2\}$}
\put(90,154){$\{1,3\}$}
\put(90,169){$\{1,4\}$}
\put(90,184){$\{2,3\}$}
\put(90,199){$\{2,4\}$}
\put(90,214){$\{3,4\}$}
\end{picture}
\caption{Two equivalent representations of the Poisson random measure $(N_t(\cdot))$: the top representation is simply the superposition of the Poisson processes $(N_t(\rho))$, $\rho \in \cP_n$; in the bottom representation, two elements of $\cL_n$ are joined at time $\sigma_i$ if they constitute the pair $R^i$. \label{fig3} }
\end{figure}

More precisely we let ${\cal{P}}_n$ denote the set of unordered pairs from 
$\cL_n = \{1, \ldots, n\}$ and for $\rho \in {\cal{P}}_n$, we write 
$\rho = \{ \rho_1, \rho_2\}$, with $\rho_1 < \rho_2$. Note that card
$(\cP_n) = n(n-1)/2$. 
Let $(N_\cdot(\rho),\ \rho \in {\cal{P}}_n)$ be independent rate one Poisson processes.
For $A \subseteq {\cal{P}}_n$ let $N_t(A) = \sum_{\rho \in A} N_t(\rho)$. 
This defines a Poisson random measure such 
that for fixed $A$, $(N_t(A),\ t \geq 0)$ is a Poisson process with intensity 
card$(A)$. Let $\sigma_1 < \sigma_2 < \cdots$ be the jump times of $(N_t({\cal{P}}_n),\ t \geq 0)$, 
and $R^i = \{R^i_1, R^i_2\}$ be the pair corresponding to 
time $\sigma_i.$ Two possible representations of this Poisson random measure are shown in Figure \ref{fig3}.

For $\ell \in \cL_n,$ let ${\cal{P}}^{(\ell)} \subseteq {\cal{P}}_n$ be the set of pairs that 
contain $\ell,$ so that card$({\cal{P}}^{(\ell)}) = n-1.$ Let 
$\tau^\ell_1 < \tau^\ell_2 < \cdots$ be the jump times of 
$(N_t({\cal{P}}^{(\ell)}),\ t \geq 0).$ We write $N_t(\ell)$ instead of 
$N_t({\cal{P}}^{(\ell)})$. Note that  
\begin{equation*}
   \sum_{\rho \in {\cal{P}}_n} N_t (\rho) = N_t(\cP_n) = \frac{1}{2}
   \sum_{\ell \in {\cal{L}}_n} N_t(\ell).
\end{equation*}

\begin{figure}[h]
\begin{picture}(400,200)(0,0)
\put(120,20){\includegraphics[scale=1,angle=0]{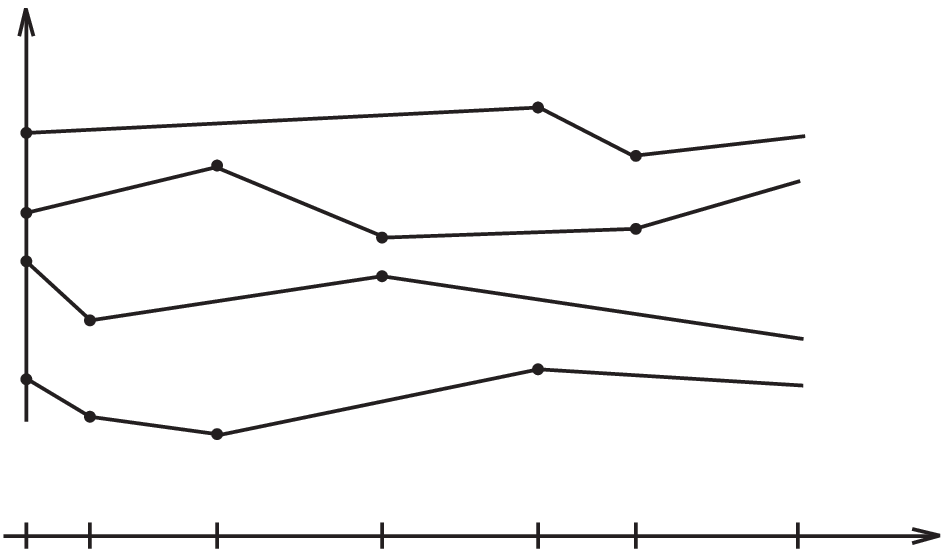}}
\put(124,25){$0$}
\put(140,25){$\tau^1_1$} \put(178,25){$\tau^1_2$}
\put(224,25){$\tau^2_2$}\put(270,25){$\tau^1_3$}
\put(298,25){$\tau^3_2$}
\put(140,2){$\tau^2_1$} \put(178,2){$\tau^3_1$}
\put(224,2){$\tau^3_2$}\put(270,2){$\tau^4_1$}
\put(298,2){$\tau^4_2$}
\put(142,13){\line(0,1){7}}\put(144,13){\line(0,1){7}}
\put(180,13){\line(0,1){7}}\put(182,13){\line(0,1){7}}
\put(226,13){\line(0,1){7}}\put(228,13){\line(0,1){7}}
\put(272,13){\line(0,1){7}}\put(274,13){\line(0,1){7}}
\put(300,13){\line(0,1){7}}\put(302,13){\line(0,1){7}}
\put(350,25){$t$}
\put(111,85){$x_1$}
\put(111,119){$x_2$}
\put(111,134){$x_3$}
\put(111,157){$x_4$}
\put(110,190){$\IR^d$}
\put(158,78){$X^1$}
\put(170,112){$X^2$}
\put(210,139){$X^3$}
\put(185,166){$X^4$}
\end{picture}
\caption{Illustration of the motion processes $X^\ell$ in the case where $n=4$ and $X^\ell_0 = x^\ell$, $\ell=1,\dots,4$. \label{fig4} }
\end{figure}

We now define the motion process needed.
For $\ell \in {\cal{L}}_n$ and $i \geq 0$, let 
$(\tilde{X}_t^{\ell, (i)},\ t \geq 0)$ be i.i.d.~copies of the process 
$(\tilde{X}_t)$ defined before Example \ref{rdexample2.1}. Set 
\begin{equation*}
   X^\ell_t = \left\{ \begin{array}{ll}
     X^\ell_0 + \tilde{X}_t^{\ell,(1)}, & 0 \leq t \leq \tau^\ell_1,\\
    \\
     X^\ell_{\tau^\ell_i} + \tilde{X}^{\ell, (i+1)}_{t-\tau^\ell_i}, 
          & \tau^\ell_i < t < \tau^\ell_{i+1}.
     \end{array}\right.
\end{equation*} 
This motion process is illustrated in Figure \ref{fig4}.

It will be useful to define $X^\ell_t$ for certain $t<0$. 
For given $(t_1,x_1), \ldots, (t_n,x_n)$, under the measure  
$P_{(t_1, x_1), \ldots, (t_n, x_n)}$ we set 
\begin{equation*}
 X^\ell_t = \tilde{X}^{\ell,(0)}_{t+t_\ell} \quad \mbox{ for $-t_\ell \leq t \leq 0$.}
\end{equation*}
Finally we set $\tau^\ell_0 = -t_\ell$.
The following theorem gives a formula for the $n$-th moments, and it is 
the main result of this section.  
\begin{thm} 
\label{nth-moment}
The $n$-th moments are given by
\begin{eqnarray}
\label{rdeq4.4}
\lefteqn{ E[u(t,x_1) \cdots u(t, x_n)]  }\\
   &=& e^{tn(n-1)/2} E_{(0,x_1), \ldots, (0, x_n)}
    \left[\prod^{N_t({\cal{P}}_n)}_{i=1}
      f(X_{\sigma_i}^{R^i_1}-X^{R^i_2}_{\sigma_i})\right.  
    \nonumber\\
   && \qquad \qquad \times \left.\prod_{\ell \in {\cal{L}}_n} \prod^{N_t(\ell)}_{i=1}
    S(\tau_i^\ell-\tau^\ell_{i-1}, \IR^d) \cdot \prod_{\ell \in {\cal{L}}_n}
     w(t-\tau^\ell_{N_t(\ell)}, X^\ell_{\tau_{N_t(\ell)}}) \right].
     \nonumber
\end{eqnarray}
\end{thm}
The proof of this theorem requires some preliminaries. Let $I_m(t,x)$, $m \geq 0$, be as defined 
in (\ref{rd3.3a}). For $0 \leq s \leq t,$ set
\begin{equation*}
I_{m+1}(s,t,x) = \int_0^s \int_{\IR^d} S(t-r, x-y) I_m (r, y) F(dr, dy),
\end{equation*}
so that $I_m(t,t,x) = I_m(t,x)$ for $ m \geq 1$. For $m=0$ and $0 \leq s < t,$ we 
set $I_0(s,t,x) = I_0(t,t,x) = w(t,x).$ Let 
\begin{equation*}
   I(s; (m_i, t_i, x_i), i = 1, \ldots, n) = E\left[\prod^n_{i=1} I_{m_i} (s, t_i, x_i)\right],
\end{equation*}
for $m_i \geq 0$, $s \leq \min(t_1, \ldots, t_n)$, $x_i \in \IR^d$, 
$i = 1, \ldots, n.$ We begin by giving an inductive expression for this 
expectation. 

\begin{lem} Suppose $m_1 + \cdots + m_n = m.$
\begin{itemize}
\item[(a)] If $m = 0,$ then 
\begin{equation*}
   I(s; (0, t_i, x_i), i = 1, \ldots, n) = \prod^n_{i=1} w(t_i, x_i).
\end{equation*}
\item[(b)] If $m \geq 1,$ then
\begin{eqnarray}\nonumber
\lefteqn{ I(s;(m_i, t_i, x_i), i=1, \ldots, n)  }\\ \nonumber
  &=& \sum_{\rho \in \cP_n:\, m_{\rho_1}\cdot m_{\rho_2} > 0}
     \int_0^s dr \int_{\IR^d} S(t_{\rho_1}-r, dy_1)\int_{\IR^d} S(t_{\rho_2}-r, dy_2)
     f(x_{\rho_1}-y_1-x_{\rho_2}+y_2) \\ \nonumber \\
  && \qquad \qquad \times \,   E\left[\prod_{i=1}^2 I_{m_{\rho_i}-1}(r,r,x_{\rho_i}-y_i) \cdot  
    \prod_{k \in \cL_n \setminus \rho} I_{m_k}(r,t_k, x_k)\right].
\label{rdeq4.0}
\end{eqnarray}
\label{rdlem4.1}
\end{itemize}
\end{lem}

\proof Part (a) follows immediately from the definitions. For part (b),
if $m =1$, then $n-1$ of the $m_i$ are equal to $0$ and so $n-1$ of the 
$I_{m_i}(s,t_i,x_i)$ are deterministic. The one $I_{m_i}(s,t_i,x_i)$ with 
$m_i =1$ is a martingale with mean zero, implying that
$I(s;(m_i, t_i, x_i), i=1, \ldots, n)=0$. The expression in formula 
(\ref{rdeq4.0}) is also equal to $0$ since there is no $\rho \in \cP_n$ such 
that $m_{\rho_1}\cdot m_{\rho_2} > 0$. 
   
If $m\geq 2$, we distinguish two cases. The first case is where all but one
of the $m_i$ are zero. In this 
case, $I(s;(m_i, t_i, x_i), i=1, \ldots, n)$ and expression (\ref{rdeq4.0}) 
vanish, for the same reasons as in the case $m=1$. We now consider the second 
case, in which there is at least one $\rho \in \cP_n$ with 
$m_{\rho_1}\cdot m_{\rho_2} > 0$. 

Using the smoothed kernels $S_\ell = \psi_\ell \ast S$, as in section
\ref{sec4.1}, we define 
\begin{equation*}
I^\ell_{m+1}(s,t,x) = \int^s_0 \int_{\IR^d} S_\ell(t-r, x-y) I_m(r,y) F(dr, dy).
\end{equation*}
For fixed $(t_i, x_i)$, $s \mapsto I^\ell_{m_i}(s, t_i, x_i)$ is a martingale, and 
according to \cite[Thm. 2.5]{walsh}, if $m_i > 0$ and $m_j > 0$, then the 
mutual variation process of $I^\ell_{m_i} (\cdot, t_i, x_i)$ and 
$I^\ell_{m_j}(\cdot, t_j, x_j)$ is 
\begin{eqnarray*}
s & \mapsto & \int_0^s dr \int_{\IR^d} dy_1\, S_\ell (t_i-r, x_i-y_1) \int_{\IR^d} dy_2\, 
S_\ell (t_j-r, x_j-y_2) \\
   && \qquad \qquad \times \, f(y_1-y_2)I_{m_i-1}(r,r,y_1)I_{m_j-1}(r,r,y_2).
\end{eqnarray*}
We now apply It\^o's formula \cite[Theorem 5.10]{CW} to the function 
$f(a_1, \ldots, a_n) = a_1 \cdots a_n$ and the $n$ martingales 
$I^\ell_{m_i}(\cdot, t_i, x_i)$, $i = 1, \ldots, n$. Note that 
\begin{equation*}
\frac{\partial^2f}{\partial a_i^2} = 0 \qquad\mbox{  and  }\qquad \frac{\partial^2f}{\partial a_i \partial a_j} 
= \prod_{k \in \{1, \ldots, n\} \setminus \{i,j\}} a_k \quad \mbox{  if  } i \not= j.
\end{equation*}
The stochastic integrals terms given by It\^o's formula have mean zero, because the 
$I_{m_i}(\cdot, \cdot, \cdot)$ have bounded moments of all orders, so 
taking expectations we reach 
\begin{eqnarray}\nonumber
   && E\left[\prod^n_{i=1} I_{m_i}^\ell(s, t_i, x_i)\right] \\
   &&\qquad = \sum_{\rho \in {\cal{P}}_n:\, m_{\rho_1}\cdot m_{\rho_2} > 0} \int_0^s dr 
\int_{\IR^d} dy_1\, S_\ell(t_{\rho_1}-r, x_{\rho_1}-y_1) \int_{\IR^d} dy_2\, 
S_\ell(t_{\rho_2}-r, x_{\rho_2}-y_2) \nonumber \\
&& \qquad\qquad \qquad \times \; f(y_1-y_2)
     E\left[I_{m_{\rho_1}-1}(r,r,y_1)I_{m_{\rho_2}-1}(r,r,y_2) 
       \prod_{k \in {\cal{L}}_n \setminus \rho} I^\ell_{m_k}(r,t_k, x_k)\right].
\label{rdeq4.2}
\end{eqnarray}
The variables $I_{m_k}(r,t_k, x_k)$ and  $I^\ell_{m_k}(r,t_k, x_k)$ are both bounded in $L^p$ for all $p$
and continuous in $L^2$ in $(r, x_k)$, so that they are 
continuous in $L^p$ in the variables $(r, x_k)$. This implies that the expectation in (\ref{rdeq4.2}) 
is continuous in $(r, x_1, \ldots, x_n)$. Using the change of variables 
$z_1 = x_{\rho_1}-y_1$ and $z_2 = x_{\rho_2}-y_2,$ we let $\ell \to \infty$ 
in (\ref{rdeq4.2}). The left-hand side converges to 
$I(s; (m_i, t_i, x_i), i=1,\ldots, n)$ and the right-hand side converges to 
formula (\ref{rdeq4.0}), completing the proof. 
\hfill $\Box$
\bigskip

Define 
\begin{eqnarray*}
 && J(t; (m_i, t_i, x_i), i = 1, \ldots, n) \\
 &&\qquad = e^{tn(n-1)/2} E_{(t_1, x_1), \ldots, (t_n, x_n)}\left[1_{\{N_t(\ell)=m_\ell, 
    \ \ell \in {\cal{L}}_n\}} \prod^{\frac{1}{2}(m_1+ \cdots + m_n)}_{i=1} f(X^{R^i_1}_{\sigma_i} 
      - X^{R^i_2}_{\sigma_i})\right. \\
   \\
   && \qquad\qquad \times \left.\prod_{\ell \in {\cal{L}}_n} \prod^{m_i}_{i=1}
 S(\tau^\ell_i - \tau^\ell_{i-1}, \IR^d) \cdot \prod_{\ell \in {\cal{L}}_n} 
w(t- \tau^\ell_{m_\ell}, X^\ell_{\tau^\ell_{m_\ell}})\right].
\end{eqnarray*}
The next aim is to show that these expectations satisfy a similar
inductive formula.
\begin{lem} Suppose $m_1 + \cdots + m_n = m.$
\begin{itemize}
\item[(a)] If $m = 0,$ then
\begin{equation*}
J(t; (0, t_\ell, x_\ell), \ell = 1, \ldots, n) = \prod^n_{\ell=1} w(t+t_\ell, x_\ell).
\end{equation*}
\item[(b)] If $m \geq 1,$ then $J(t; (m_\ell, t_\ell, x_\ell), \ell \in {\cal{L}}_n)$ is equal to
\begin{eqnarray} \nonumber
  && \sum_{\rho \in {\cal{P}}_n:\, m_{\rho_1}\cdot m_{\rho_2} > 0} \int^t_0 ds \int_{\IR^d} 
S(t_{\rho_1} + s, dy_1) \int_{\IR^d} S(t_{\rho_2} + s, dy_2) 
f(x_{\rho_1} - y_1 - x_{\rho_2} + y_2)\\ \nonumber
\\
  &&\qquad \qquad \times \; J(t-s; (m_{\rho_i}-1, 0, x_{\rho_i}-y_i),\ i=1,2; 
(m_\ell, s+t_\ell, x_\ell),\ \ell \in {\cal{L}}_n \setminus \rho). 
\label{rdeq4.3}
\end{eqnarray} \label{rdlem4.2}
\end{itemize}
\end{lem}

\proof Part (a) follows immediately from the definitions. For part (b), 
in the case that only one of the $m_i$ are non-zero then 
$J(t; (m_\ell, t_\ell, x_\ell), \ell \in {\cal{L}}_n) = 0$  since 
$P_{(t_1,x_1),\dots,(t_n,x_n)} \{N_t(\ell) = m_\ell,\ \ell \in {\cal{L}}_n\} =0$,
and formula (\ref{rdeq4.3}) is also equal to $0$ since there is no 
$\rho \in \cP_n$ with $m_{\rho_1}\cdot m_{\rho_2} > 0$.
   
We now suppose that $m \geq 2$ and that there is at least one 
$\rho \in \cP_n$ with $m_{\rho_1}\cdot m_{\rho_2} > 0$. 
In this case $\{N_t(\ell) = m_\ell,\ \ell \in {\cal{L}}_n\} \subset \{\sigma_1 \leq t\}$, 
and we are going to use the Markov property of $N_t({\cal{P}}_n)$ at 
time $\sigma_1$. Indeed, 
choosing $\FC_1= \sigma \{ \sigma_1, R^1, X^{R^1_1}_{\sigma_1}, X^{R^1_2}_{\sigma_1}\}$,  
we may rewrite $J(t; (m_\rho, t_\ell, x_\ell),\ \ell \in {\cal{L}}_n)$ as 
\begin{eqnarray*} 
 && \sum_{\rho \in {\cal{P}}_n:\, m_{\rho_1}\cdot m_{\rho_2} > 0} E_{(t_1, x_1),
        \ldots, (t_n, x_n)} \left[ 1_{\{\sigma_1 \leq t, R^1 = \rho\}}
         f(X^{\rho^1}_{\sigma_1}-X^{\rho^2}_{\sigma_1}) 
         e^{t n(n-1)/2} \prod_{\ell \in \rho} S(\tau_1^\ell-\tau_0^\ell, \IR^d) \right. \\
   && \qquad \times \, E_{(t_1, x_1), \ldots, (t_n, x_n)} [1_{\{N_t(\ell) - N_{\sigma_1}(\ell)
       = m_\ell, \ \ell \in {\cal{L}}_n \setminus \rho\} \cap \{N_t(\ell) -
        N_{\sigma_1}(\ell) =
        m_\ell-1,\ \ell \in \rho\}} \prod^m_{\ell = 2} f(X^{R^i_1}_{\sigma_i} -
        X^{R^i_2}_{\sigma_i})  \\
   && \qquad \times \left. \prod_{\ell \in {\cal{L}}_n \setminus \rho}\ \prod^{m_i}_{i=1}
        S(\tau^\ell_i - \tau^\ell_{i-1}, \IR^d) \cdot \prod_{\ell \in \rho}
        \prod^{m_i}_{i=2} S(\tau^\ell_i - \tau^\ell_{i-1}, \IR^d) \cdot \prod_{\ell \in
         {\cal{L}}_n} w(t-\tau^\ell_{m_\ell}, X^\ell_{\tau^\ell_{m_\ell}} ) \vert
          \FC_{1}]\right].
\end{eqnarray*}
Note that at time $\sigma_1,$ on $\{R^1 = \rho\}$, the processes $X^{\rho_i}$ start 
afresh from $X^{\rho_i}_{\sigma_1}$, $i=1,2$, while for 
$\ell \in {\cal{L}}_n \setminus \rho$, $X^\ell$ has seen no jump from $-t_\ell$ 
to $\sigma_1$, that is for $\sigma_1+t_\ell$ units of time. 
Using the strong Markov property at $\sigma_1$, the conditional expectation above 
multiplied by $e^{(t-\sigma_1)n(n-1)/2}$ is equal to 
\begin{equation*}
   J(t-\sigma_1; (m_{\rho_i}-1, 0, X^{\rho_i}_{\sigma_1}),\ i=1,2; 
   (m_\ell, \sigma_1 + t_\ell, x_\ell),\ \ell \in {\cal{L}}_n \setminus \rho).
\end{equation*}
Therefore, $J(t; (m_\ell, t_\ell, x_\ell), \ell \in {\cal{L}}_n) $ is equal to
\begin{eqnarray*}
   && \sum_{\rho \in {\cal{P}}_n:\, m_{\rho_1}\cdot m_{\rho_2} > 0} 
       E_{(t_1, x_1), \ldots, (t_n, x_n)} [ 1_{\{\sigma_1 \leq t, R^1 = \rho\}}
       e^{\sigma_1 n(n-1)/2} f(X^{\rho_1}_{\sigma_1} - X^{\rho_2}_{\sigma_1}) \prod_{\ell \in \rho}
       S(\tau_1^\ell - \tau^1_0, \IR^d)\\
&& \qquad\qquad\qquad \times   J(t-\sigma_1; (m_{\rho_i}-1, 0, X^{\rho_i}_{\sigma_1}),\
      i=1,2; (m_\ell, \sigma_1 + t_\ell, x_\ell),\ 
      \ell \in {\cal{L}}_n \setminus \rho)].
\end{eqnarray*}
The variable $\sigma_1$ is exponential with mean $2/(n(n-1))$ and the variable $R$ 
is independent and uniformly distributed over ${\cal{P}}_n$.
Taking the expectation over $\sigma_1, R^1, X^{R^1_i}_{\sigma_1}$ we reach (\ref{rdeq4.3}).
\hfill $\Box$ 
\bigskip

\pf {\it of Theorem \ref{nth-moment}.} We note that it suffices to prove, when $t_i \geq t$, that
\begin{equation}\label{rdeq4.5}
I(t; (m_1, t_1, x_1), \ldots, (m_n, t_n, x_n)) = J(t; (m_1, t_1-t, x_1), \ldots, (m_n, t_n-t, x_n)).
\end{equation}
Indeed, in this case, by Proposition \ref{prop3.1},
\begin{eqnarray*}
   E[u(t, x_1) \cdots u(t, x_n)] &=& \sum^\infty_{m_1=0} \cdots \sum^\infty_{m_n=0}
      E(I_{m_1}(t, t, x_1) \cdots I_{m_n}(t, t, x_n))\\
   \\
   &=& \sum^\infty_{m_1=0} \cdots \sum^\infty_{m_n=0} I(t; (m_1, t, x_1), \ldots, 
      (m_n, t, x_n))\\
   \\
   &=& \sum^\infty_{m_1=0} \cdots \sum^\infty_{m_n=0} J(t; (m_1, 0, x_1), \ldots, 
      (m_n, 0, x_n)),
\end{eqnarray*}
and this is equal to the expression in (\ref{rdeq4.4}).

Let $m = m_1 + \cdots + m_n$. We are going to prove (\ref{rdeq4.5}) by 
induction on $m$. If $m=0,$ then (\ref{rdeq4.5}) follows from Lemma 
\ref{rdlem4.1}(a) and Lemma \ref{rdlem4.2}(a), since both sides of 
(\ref{rdeq4.5}) are equal to $ w(t_1, x_1) \cdots w(t_n, x_n)$.
Now assume inductively that (\ref{rdeq4.5}) holds for $m-1\geq 0$. By Lemma \ref{rdlem4.2}(b), 
\begin{eqnarray}\nonumber
   && J(t; (m_\ell, t_\ell - t, x_\ell),\ \ell \in {\cal{L}}_n) \\
   &&\qquad = \sum_{\rho \in {\cal{P}}_n:\, m_{\rho_1}\cdot m_{\rho_2} > 0} 
        \int_0^t ds \int_{\IR^d} S(t_{\rho_1} - t+s, dy_1) \int_{\IR^d} 
        S(t_{\rho_2}-t+s, dy_2) \nonumber\\
   && \qquad\qquad \qquad\qquad \times \; f(x_{\rho_1}-y_1-x_{\rho_2}+y_2) \label{rdeq4.6}  \\
   &&\qquad\qquad \qquad\qquad \times \; J(t-s; (m_{\rho_i}-1, 0, x_{\rho_i}-y_i),\ i = 1,2;
       (m_\ell, s+t_\ell-t, x_\ell),\ \ell \in {\cal{L}}_n - \rho).
\nonumber
\end{eqnarray}
By the induction hypothesis, the last factor $J(t-s; \ldots)$ is equal to
\begin{equation*}
   I(t-s; (m_{\rho_i}-1, t-s, x_{\rho_i}-y_i),\ i=1,2; 
(m_\ell, t_\ell, x_\ell),\ \ell \in {\cal{L}}_n \setminus \rho).
\end{equation*}
Now use the change of variables $r=t-s$ and Lemma \ref{rdlem4.1}(b) to see 
that the right-hand side of (\ref{rdeq4.6}) is equal to 
$I(t; (m_\ell, t_\ell, x_\ell),\ \ell \in {\cal{L}}_n).$ This completes the 
proof.  
\hfill $\Box$
\bigskip

\begin{rem} The intuition behind equality (\ref{rdeq4.5}) is the following. Suppose $n=4$ and consider space-time positions $(t_1,x_1),\dots,(t_4,x_4)$, as in Figure \ref{fig5}.
\begin{figure}[h]
\begin{picture}(400,180)(0,0)
\put(120,20){\includegraphics[scale=1,angle=0]{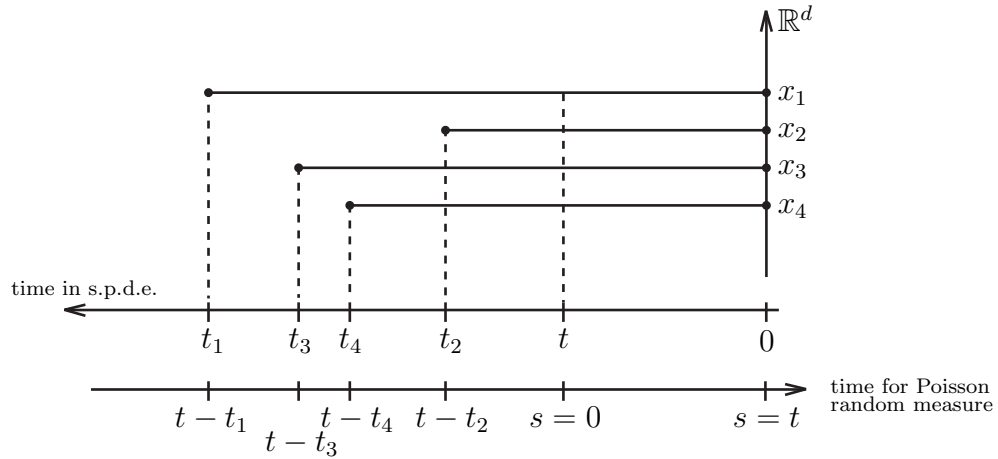}}
\put(172,40){$t_1$}
\put(262,40){$t_2$}
\put(205,40){$t_3$}
\put(224,40){$t_4$}
\put(307,40){$t$}
\put(383,39){$0$}
\put(162,10){$t-t_1$}
\put(196,01){$t-t_3$}
\put(217,10){$t-t_4$}
\put(253,10){$t-t_2$}
\put(296,10){$s=0$}
\put(373,10){$s=t$}
\put(390,91){$x_4$}
\put(390,106){$x_3$}
\put(390,120){$x_2$}
\put(390,134){$x_1$}
\put(390,160){$\IR^d$}
\put(410,23){{\scriptsize time for Poisson}}
\put(410,16){{\scriptsize random measure}}
\put(100,60){{\scriptsize time in s.p.d.e.}}

\end{picture}
\caption{Illustration for equality (\ref{rdeq4.5}), with the direction of time for the s.p.d.e.~and for the Poisson random measure. \label{fig5}}
\end{figure}
The quantity $I(t; (m_\ell,t_\ell, x_\ell),\ \ell =1,\dots 4)$ is the expected product of iterated integrals, where the left-most integral is up to time $t \leq \min(t_1,\dots,t_4)$ and the order of the iterated integrals are $m_1,\dots,m_4$.

   On the other hand, time $s$ for the Poisson random measure runs in the opposite direction as in the s.p.d.e. (see Figure \ref{fig5}). In the quantity $J(t, (m_\ell, t_\ell - t, x_\ell),\ \ell =1,\dots 4)$, the process $X^\ell$ starts at negative time $t - t_\ell$, and there are no Poisson pairs during negative time. During the time interval $s=0$ to $s=t$, the number of Poisson pairs containing $x_\ell$ is set to $m_\ell$. With these constraints, $I(t; (m_\ell,t_\ell, x_\ell),\ \ell =1,\dots 4) = J(t, (m_\ell, t_\ell - t, x_\ell),\ \ell =1,\dots 4)$ as stated in (\ref{rdeq4.5}).
\end{rem}


\begin{thebibliography}{[XX]}

\bibitem{albevario1} Albeverio, S.A. \& Hoegh-Krohn, R.J. Mathematical theory of Feynman path integrals. Lecture Notes in Mathematics, Vol. 523. Springer-Verlag, Berlin-New York, 1976. 

\bibitem{albevario2}  Albeverio, S., Blanchard, Ph., Combe, Ph., Hoegh-Krohn, R. \& Sirugue, M. Local relativistic invariant flows for quantum fields.  Comm. Math. Phys.  90  (1983), 329-351.  

\bibitem{carmmolch} Carmona, R.A. \& Molchanov, S.A. Parabolic Anderson problem and intermittency.  Mem. Amer. Math. Soc.  108  (1994),  no. 518.

\bibitem{CW}  Chung, K.L. \& Williams, R.J. Introduction to stochastic integration. Second edition. Probability and its Applications. Birkhäuser Boston, Inc., Boston, MA, 1990.  

\bibitem{dapratoz} Da Prato, G. \& Zabczyk, J. Stochastic equations in infinite dimensions. Encyclopedia of Mathematics and its Applications, 44. Cambridge University Press, Cambridge, 1992.

\bibitem{D99} Dalang, R.C. Extending the martingale measure stochastic integral with applications to spatially homogeneous s.p.d.e.'s.  Electron. J. Probab.  4  (1999), no. 6, 29 pp.

\bibitem{D99b} Dalang, R.C. Corrections to: ``Extending the martingale measure stochastic integral with applications to spatially homogeneous s.p.d.e.'s" Electron J. Probab. 4 (1999), no. 6, 29 pp.

\bibitem{dalangfrangos} Dalang, R.C. \& Frangos, N.E.
The stochastic wave equation in two spatial dimensions. Ann. Probab. 26 (1998), no. 1, 187-212.

\bibitem{dleveque} Dalang, R.C. \& L\'ev\`eque, O. Second-order hyperbolic s.p.d.e.'s driven by homogeneous Gaussian noise on a hyperplane. Transactions of the American Mathematical Society (2005, to appear).

\bibitem{DalMuel} Dalang, R.C. \& Mueller, C. Intermittency properties in a hyperbolic Anderson problem (in preparation).

\bibitem{fmeleard} Fournier, N. \&  M\'el\'eard, Sylvie. A stochastic 
particle numerical method for $3D$ Boltzmann equations without cutoff.  
Math. Comp.  71  (2002),  no. 238, 583-604.

\bibitem{Hersch} Hersch, R. Random evolutions: a survey of results and problems. Rocky Mountain J. Math. 4 (1974), 443-477.

\bibitem{Kac1} Kac, M. Some stochastic problems in physics and mathematics. Magnolia Petroleum Co. Lectures in Pure and Applied Science 2 (1956).

\bibitem{Kac2} Kac, M.A. A stochastic model related to the telegrapher's equation. Rocky Mountain J. Math.  4  (1974), 497--509.

\bibitem{KS} Karatzas, I. \& Shreve, S.E. Brownian motion and stochastic calculus. Second edition. Graduate Texts in Mathematics, 113. Springer-Verlag, New York, 1991.

\bibitem{leveque} L\'ev\`eque, O. Hyperbolic SPDE's driven by boundary
noises, PhD Thesis No. 2452, Ecole Polytechnique F\'ed\'erale de Lausanne, Switzerland, 2001.

\bibitem{mueller} Mueller, C. Long time existence for the wave equation with a noise term.  Ann. Probab.  25  (1997),  no. 1, 133-151.

\bibitem{oksendall} Oksendal, B., Vage, G. \& Zhao, H.Z. Asymptotic properties of the solutions to stochastic KPP equations.  Proc. Roy. Soc. Edinburgh Sect. A  130  (2000),  no. 6, 1363-1381.

\bibitem{pinsky} Pinsky, M.A. Lectures on random evolution. World scientific, 1991.

\bibitem{schwartz} Schwartz, L. Th\'eorie des distributions. Hermann, Paris,
1966.

\bibitem{walsh} Walsh, J.B. An introduction to stochastic partial differential equations.  Ecole d'\'et\'e de probabilit\'es de Saint-Flour, XIV---1984, Lecture Notes in Math., 1180, Springer, Berlin, 1986,  265-439.

\bibitem{WZ} Wheeden, R.L. \& Zygmund, A. Measure and integral. An introduction to real analysis. Pure and Applied Mathematics, Vol. 43. Marcel Dekker, Inc., New York-Basel, 1977.

\bibitem{wolpert} Wolpert, R.L. Local time and a particle picture for Euclidean field theory. J.~Funct.~Anal. 30-3 (1978), 341-357.
\end{thebibliography}
\bibliographystyle{alpha}

\end{document}